\documentclass{amsart}

\usepackage{amssymb}
\usepackage{amsmath}
\usepackage{amsthm}

\theoremstyle{plain}
\newtheorem{theorem}{Theorem}[section]

\newtheorem{definition}[theorem]{Definition}

\theoremstyle{definition}

\numberwithin{equation}{section}

\def\<#1>{\left\langle{#1}\right\rangle}
\newcommand{\set}[1]{{\left\{{\def\st{\;:\;}#1}\right\}}}
\newcommand{\abs}[1]{\left|{#1}\right|}
\newcommand{\Z}{\mathbb{Z}}
\newcommand{\Q}{\mathbb{Q}}
\newcommand{\C}{\mathbb{C}}
\DeclareMathOperator{\norm}{\mathcal{N}}
\DeclareMathOperator{\Tr}{\mathsf{Tr}}
\newcommand{\kro}[2]{\left(\frac{#1}{#2}\right)}
\newcommand{\tf}[2]{{^{#1}\!\!\!\;/\!_{#2}}}
\DeclareMathOperator{\thetaop}{\Theta}

\newcommand{\wl}{\omega}

\begin{document}
\title[Computation of central value]
      {Computation of central value of
       quadratic twists of modular $L-$functions}
\author[Z. Mao \and
        F. Rodriguez-Villegas \and
        G. Tornar\'\i a]
        {Zhengyu Mao \and
        Fernando Rodriguez-Villegas \and
        Gonzalo Tornar\'\i a}
\address{Department of Mathematics and Computer Science\\
         Rutgers university\\
         Newark, NJ 07102-1811}
\email{zmao@andromeda.rutgers.edu}
\address{Department of Mathematics\\
         University of Texas at Austin\\
         Austin, TX 78712}
\email{villegas@math.utexas.edu}
\email{tornaria@math.utexas.edu}

\date{}
\maketitle

\section{Introduction}

Let $f\in S_2(p)$ be a newform of weight two, prime level $p$.  If
$f(z)=\sum_{m=1}^{\infty}a(m) q^m$, where $q=e^{2\pi iz}$, and $D$ is a
fundamental discriminant, we define the twisted $L$-function
$$L(f,D,s)=\sum_{m=1}^{\infty} \frac{a(m)}{m^s}\left(\frac{D}{m}\right).$$
It will be
convenient to also allow $D=1$ as a fundamental discriminant, in which
case we write simply $L(f,s)$ for $L(f,1,s)$.

In this paper we consider the question of computing the twisted
central values $\set{L(f,D,1)\st \abs{D}\leq x}$ for some $x$.

It is well known that the fact that $f$ is an eigenform for the
Fricke involution yields a rapidly convergent series for
$L(f,D,1)$. Computing $L(f,D,1)$ by means of this series, which we
call the \emph{standard method}, takes time very roughly
proportional to $\abs{D}$ and therefore time very roughly
proportional to $x^2$ to compute $L(f,D,1)$ for $\abs{D}\leq
x$. We will see that this can be improved to $x^{\tf{3}{2}}$ by
using an explicit version of Waldspurger's theorem~\cite{wa1};
this theorem relates the central values $L(f,D,1)$ to the
$\abs{D}$-th Fourier coefficient of weight $\tf{3}{2}$ modular
forms in Shimura correspondence with $f$.

Concretely, the formulas we use  have the basic form
\begin{equation}\label{main-fmla}
  L(f,D,1)=\star\,\kappa_\mp \frac{\abs{c_\mp(\abs{D})}^2}{\sqrt{\abs{D}}},
        \qquad \mathop{\rm sign}(D)=\mp,
\end{equation}
where $\star=1$ if $p\nmid D$, $\star=2$ if $p\mid D$, $\kappa_-$
and $\kappa_+$ are positive constants independent of $D$, and
$c_-(\abs{D})$ (resp. $c_+(\abs{D})$) is the
$\abs{D}$-th Fourier coefficient of a certain modular form
$g_-$ (resp. $g_+$) of weight $\tf{3}{2}$.

Gross \cite{g} proves such a formula, and gives an explicit
construction of the corresponding form $g_-$, in the case that
$L(f,1)\not=0$.
The purpose of this paper is to extend Gross's work to all cases.
Specifically, we give an explicit construction of both $g_-$ and
$g_+$, regardless of the value of $L(f,1)$, together with the
corresponding values of $\kappa_-$ and $\kappa_+$ in \eqref{main-fmla}.
The proof of the validity of this construction will be given in a
later publication and relies partly in the results of \cite{bm}.

The construction gives $g_-$ and $g_+$ as linear combinations of
(generalized) theta series associated to positive definite ternary
quadratic forms. Computing the  Fourier coefficients of these theta
series up to $x$ is tantamount to running over all lattice points in
ellipsoid of volume proportional to $x^{\tf{3}{2}}$. Doing this takes time
roughly proportional to $x^{\tf{3}{2}}$ which yields our claim above.

This approach to computing $L(f,D,1)$ has several other advantages
over the standard method. First, the numbers $c(\abs{D})$ are
algebraic integers and are computed with exact arithmetic. Once
$c(\abs{D})$ is know it is trivial to compute $L(f,D,1)$ to
any desired precision. Second, the $c(\abs{D})$'s have extra
information; if $f$ has coefficients in $\Z$, for example,
\eqref{main-fmla} gives a specific square root of $L(f,D,1)$ (if
non-zero), whose sign remains a mystery.

Moreover, the actual running time of our method vs. the standard
method is, in practice, significantly better even for small $x$.


\section{Quaternion algebras and Brandt matrices}

A quaternion algebra $B$ over a field $K$ is a central simple algebra
of dimension $4$ over $K$. When $2\neq 0$ in $K$ we can give $B$
concretely by specifying a $K$-basis $\set{1,i,j,k}$ such that
\[
    i^2=\alpha,\qquad j^2=\beta,\qquad\text{and}\qquad k=ij=-ji,
\]
for some non-zero $\alpha,\beta \in K$. If $K=\Q$ we typically rescale
and assume that $\alpha, \beta \in\Z$.  A general element of $B$
then has the form $b = b_0 + b_1 i + b_2 j + b_3 k$, with $b_i\in K$
and multiplication in $B$ is determined by the above defining
relations and $K$-linearity.

 The {\it conjugate} of $b$ is defined as
\[
    \overline{b} = b_0 - b_1 i - b_2 j - b_3 k.
\]
We define  the (reduced) norm and trace of $b$ by
\[
   \norm b := b \overline{b}
   = b_0^2 - \alpha b_1^2 - \beta b_2^2 + \alpha\beta b_3^2
   , \qquad \Tr b := b + \overline{b} = 2 b_0.
\]

Let $B$ be a quaternion algebra over $K=\Q$.  For $\nu$ a rational
prime we let $\Q_\nu$ be the field of $\nu$-adic numbers and for
$\nu=\infty$ we let $\Q_\nu=\mathbb R$.  We call $\nu$, a rational
prime or $\infty$, a {\it place} of $\Q$.

The localization $B_\nu := B\otimes\Q_\nu$ is a quaternion algebra
over $\Q_\nu$.  It is a fundamental fact of Number Theory that
$B_\nu$ is either isomorphic to the algebra $M_2(\Q_\nu)$ of
$2\times2$ matrices, or a division algebra, which is unique up to
isomorphism.  (A {\it division algebra} is an algebra in which every
non-zero element has a multiplicative inverse.) The two options are
encoded in the {\it Hilbert symbol} $(\alpha,\beta)_\nu$, defined as
$+1$ if $B_\nu$ is a matrix algebra, $-1$ if it is a division algebra.
In the first case we say that $B$ is \emph{split} at $\nu$, in the
second, that $B$ is \emph{ramified} at $\nu$.

For example, if $\nu=\infty$ so $\Q_\nu = \mathbb R$ then
$(\alpha,\beta)_\infty=-1$ if and only if $\alpha<0,\beta<0$ in which
case $B_\infty$ is isomorphic to the usual Hamilton quaternions.  A
quaternion algebra $B$ is {\it definite} if it ramifies at $\infty$
otherwise it is {\it indefinite} (this notation is consistent with the
nature of the quadratic form on $B_\infty$ determined by the norm
$\norm$).

A quaternion algebra $B$ is ramified at a finite number of places and
the total number of ramified places must be even (e.g. the {\it
  Hilbert reciprocity law} says that $\prod_\nu(\alpha,\beta)_\nu=1$).
The set of ramified places determines $B$ up to isomorphism (the {\it
  local-global principle}). For any finite set $S$ with an even number
of places there is a (unique up to isomorphism) $B$ which ramifies
exactly at places in $S$.

Let $B$ be a quaternion algebra over $\Q$.  An {\it order} in $B$
is a (full rank) lattice $R\subseteq B$ which is also a ring with
$1\in R$.  As for number fields, an element of an order must be
integral over $\Z$, i.e., must satisfy a monic equation with
coefficients in $\Z$ (or even more concretely must have integral
trace and norm).  Unlike in the commutative case, however, the set of
all integral elements of $B$ is \emph{not} a ring. The best next thing
is to consider maximal orders (which always exist), i.e., orders not
properly contained in another order. But maximal orders are not
unique. In fact, if $B$ is definite, a maximal order is in general not
even unique up to isomorphism though there always is only a finite
number of isomorphism classes of maximal orders in $B$.

As an illustration consider the classical case $\alpha=\beta=-1$ of
the Hamilton quaternions.  The algebra is definite and hence ramifies
at $\nu=\infty$. It must ramify at a least one other prime, which
turns out to be only $\nu=2$. To see this note that

\[
\norm(b_0 + b_1 i + b_2 j + b_3 k)=b_0^2+b_1^2+b_2^2+b_3^2.
\]
There always is a non-trivial solution to the congruence
$b_0^2+b_1^2+b_2^2+b_3^2\equiv 0 \bmod p$ for $p$ prime. If $p$ is odd
we can lift this solution to a solution in $\Z_p$ by Hensel's
lemma obtaining a non-zero quaternion in $B_p$ of zero norm. This
implies that $B_p$ cannot be a divison algebra and hence $(-1,-1)_p=1$
for $p$ odd.  We must necessarily have then that $(-1,-1)_2=-1$.

If we want to study the representation of numbers as sum of four
squares it is natural to consider, as Lipschitz did, the arithmetic of
the quaternions with $b_i\in \Z$. These quaternions form an order
$R'$, but, as it turns out, it is not maximal. Indeed, as Hurwitz
noted, $\rho:=\tfrac12(1+i+j+k)$ is integral ($\norm \rho=1$ and
$\Tr \rho =1$) and $R:=R'+\Z\rho$ is also an order of $B$
strictly containing $R'$.

Moreover, $R$ is maximal and hence its arithmetic is significantly
simpler than that of $R'$. Hurwitz showed, for example, that there is
a left and right division algorithm in $R$, from which it follows that
every positive integer is a sum of four squares.

Fix a prime $p$ and let $B$ be the quaternion algebra over $\Q$
ramified precisely at $\infty$ and $p$. Let $R$ be a fixed maximal
order in $B$.  A right ideal $I$ of $R$ is a lattice in $B$ that is
stable under right multiplication by $R$. Two right ideals $I$ and $J$
are in the same {\it class} if $J=bI$ with $b\in B^\times$.  The set
of right ideal classes is finite; let $n$ be its number. Chose a set
of representatives $\set{I_1,\dotsc,I_n}$ of the classes. (We should
emphasize here that contrary to the commutative setting there is {\it no}
natural group structure on the set of classes.)

Consider the vector space $V$ of formal linear combinations
\[
\sum_{i=1}^n a_i\,[I_i], \qquad a_i \in \C
\]
(here $[I]$ denotes the class of $I$).

For each integer $m$ there is an $n\times n$ matrix $B_m$ acting on
$V$. Pizer~\cite{p} gives an efficient algorithm for computing these
{\it Brandt matrices}: its coefficients are given by the
representation numbers of the norm form for certain quaternary
lattices in $B$.

The Brandt matrices commute with each other and are self-adjoint with
respect to the \emph{height pairing} on $V$ (see \S1 and \S2
of~\cite{g} for an account of this.) From this it follows that there
is basis of $V$ consisting of simultaneous eigenvectors of all $B_m$.

It follows from Eichler's trace formula that there is a one to one
correspondence between Hecke eigenforms of weight $2$ and level $p$
(cf.  \cite[\S5]{g}) and eigenvectors in $V$ of all Brandt matrices
(up to a constant multiple).

If $f$ is the Hecke eigenform we let $e_f$ be the corresponding
eigenvector (well defined up to a constant). Then $B_me_f=a_me_f$
where $T_mf =a_mf$ and $T_m$ is the $m$-th Hecke operator.

\section{Construction of $g_-$ and $g_+$}

Let $e_f$ be the eigenvector for the Brandt matrices for $R$
corresponding to $f$ as in the last section.
One can use linear algebra to find its coefficients
\[
    e_f=\sum_{i=1}^n a_i [I_i],
\]
by computing the Brandt matrices, and
from the knowledge of a few
eigenvalues (i.e.  Fourier coefficients) of $f$.

We will describe below the construction of certain generalized theta
series $\thetaop_{l^\ast}([I_i])$ corresponding to each ideal class
$[I_i]$, and then define
\[
   \thetaop_{l^\ast}(e_f)
      := \sum_{i=1}^n a_i\,\thetaop_{l^\ast}([I_i])
       = \sum_{n=1}^{\infty} c_{l^\ast}(n)\, q^n.
\]
Here $l^\ast$ is a fundamental discriminant for which we will consider
three cases: $l^\ast=1$, which is Gross's construction of $g_-$;
$l^\ast=l$ for an odd prime $l\neq p$ such that $l\equiv 1\pmod{4}$,
which will generalize Gross's construction of $g_-$; and $l^\ast=-l$
for an odd prime $l\neq p$ such that $l\equiv 3\pmod{4}$, which will
give a construction of $g_+$.

Furthermore, for any fundamental discriminant $D$ such that
$Dl^\ast<0$, the following formula holds
\begin{equation}\label{eq:formula}
   L(f,l^\ast,1)\,L(f,D,1)
    = \star\,\kappa_f
            \frac{\abs{c_{l^\ast}(\abs{D})}^2}
                 {\sqrt{\abs{Dl^\ast}}},
\end{equation}
where $\star=1$ if $p\nmid D$, $\star=2$ if $p\mid D$,
and $\kappa_f := \frac{\<f,f>} {\<e_f,e_f>}$
is a positive constant independent of $D$ or $l^\ast$. Here
$\<e_f,e_f>$ is the height of $e_f$, and $\<f,f>$ is the
Petersson norm of $f$ (cf. \S4 and \S7 of~\cite{g}.)
For $l=1$, it was proved by Gross in~\cite[Proposition 13.5]{g}.
The proof of this formula for the case $l\neq 1$ will be given in a
later publication.

Note that, as a corollary, we have $\thetaop_{l^\ast}(e_f)\neq 0$ if and
only if $L(f,l^\ast,1)\neq 0$, and this happens for infinitely many
$l^\ast>0$ and for infinitely many $l^\ast<0$, as follows
from~\cite{bfh}.


\subsection{Gross's construction of $\thetaop_1$}

Let $R_i:=\set{b\in B\st bI_i\subset I_i}$ be the left order of
$I_i$. The $R_i$ are maximal orders in $B$, and each conjugacy
class of maximal orders has a representative $R_i$ for some $i$.

We let $S_i^0:=\set{b\in\Z+2R_i\st\Tr b = 0}$, a ternary
lattice,
and define
\[
  \thetaop_1([I_i])
    := \frac{1}{2}\sum_{b\in S_i^0} q^{\,\norm b}.
\]
Then $\thetaop_1([I_i])$ is a weight $\tf{3}{2}$ modular form
of level $4p$ and trivial character.


\subsection{Weight functions and $\thetaop_l$}


Fix an odd prime $l\neq p$. In order to generalize Gross's
method, we need to construct certain weight functions
$\wl_l(I_i,\cdot)$ on $S_i^0$ with values in $\set{0,\pm 1}$.
There is a choice of sign in the construction, and some care is
needed to ensure that the choice is consistent from one ideal to
another.
It will be the case that $\wl_l(I_i,b)=0$ unless $l\mid\norm b$,
and thus we define a \emph{generalized theta series}
\[
  \thetaop_l([I_i])
    := \frac{1}{2}\sum_{b\in S_i^0} \wl_l(I_i,b) q^{\,\norm b/l},
\]
a modular form of weight $\tf{3}{2}$ and level $4p$ with trivial
character. In addition, $\thetaop_l([I_i])$ is already a cusp form
whenever $l\neq 1$, although it might be zero.

\begin{definition}
Given a pair $(L,v)$, where $L$ is an integral $\Z_l$-lattice
of rank $3$ with $l\nmid\det L$, and $v\in L$ is such that
$l\mid\norm v$ but $v\not\in l\,L$, we define its \emph{weight
function} $\wl_{l,v}: L \rightarrow \set{0,\pm 1}$ to be
\[
 \wl_{l,v}(v') := \begin{cases}
   0 & \text{ if $l\nmid\norm v'$,} \\
   \chi_l(\<v,v'>) & \text{ if $l\nmid\<v,v'>$,} \\
   \chi_l(k) & \text{ if $v'-k\,v\in l\,L$.}
  \end{cases}
\]
Here $\chi_l$ is the quadratic character of conductor $l$,
and $\norm v := \frac{1}{2}\<v,v>$.
\end{definition}

This is well defined, because if $v,v'\in L$ are such that $\norm
v\equiv\norm v'\equiv\<v,v'>\equiv 0\pmod{l}$, then $v$ and $v'$
must be collinear modulo $l$, since $L$ is unimodular. This means
that, assuming $v\not\in l\,L$, there is indeed a well defined
$k\in\Z/l\,\Z$ such that $v'-k\,v\in l\,L$.

Note that there are, for different choices of $v$, two different
weight functions for each $L$, opposite to each other; the definition
above singles out the one for which $\wl_{l,v}(v)=+1$.

We will apply the above definition to the ternary lattices
$S_i^0(\Z_l):=S_i^0\otimes\Z_l$. Fix a quaternion $b_0\in
S^0:=\set{b\in\Z+2R\st\Tr b = 0}$, and such that
$l\mid\norm b_0$ but $b_0\not\in lS^0$. For each $I_i$, find
$x_i\in I_i$ such that $l\nmid n_i:=\norm x_i/\norm I_i$.  Then
$x_i$ is a local generator of $I_i$, and
$b_i:=x_i\,b_0\,x_i^{-1}\in S_i^0(\Z_l)$.
We finally set
\[
   \wl_l(I_i, b) := \chi_l(n_i)\,
      \wl_{l,b_i}(b),
\]
where $\wl_{l,b_i}$ is the weight function of the pair
$(S_i^0(\Z_l),b_i)$.

\subsection{Odd weight functions and $\thetaop_{-l}$}

When $l\equiv 3\pmod{4}$ the weight functions $\wl_l(I_i,\cdot)$ are
odd, since $\chi_l$ is odd. Therefore, we will have $\thetaop_l=0$.
To address this problem, we will construct a different kind of
weight functions $\wl_p(I_i,\cdot)$, and then define
\[
  \thetaop_{-l}([I_i])
    := \frac{1}{2}\sum_{b\in S_i^0} \wl_p(I_i,b)\,\wl_l(I_i,b) q^{\,\norm b/l},
\]
which will be a modular form of weight $\tf{3}{2}$, this time of level $4p^2$.
Again, $\thetaop_{-l}([I_i])$ is a cusp form, which might be zero.
Note that we could have used the product of two odd weight functions
$\wl_{l_1}$ and $\wl_{l_2}$, but this construction would only lead us to
the same $g_-$. By using the weight functions $\wl_p$ we get a
construction of $g_+$ instead.

\begin{definition}
Given a triple $(L,v,\psi)$ where $L$ is an integral
$\Z_p$-lattice of rank $3$ with level $p$ and determinant
$p^2$ (i.e. $L$ is $\Z_p$-equivalent to $S^0(\Z_p)$,)
$v\in L$ is such that $p\nmid\norm v$, and
$\psi:\Z/p\,\Z\rightarrow\C$ is a periodic function
modulo $p$, the \emph{weight function}
$\wl_{\psi,v}:L\rightarrow\C$ is defined by
\[
  \wl_{\psi,v}(v') := \psi(\<v,v'>).
\]
\end{definition}

Clearly, the weight function $\wl_{\psi,v}$ will be odd if and
only if $\psi$ itself is odd. From now on we assume that $\psi$
is a fixed odd periodic function such that
\begin{equation}\label{eq:psimod1}
\abs{\psi(t)}=1\qquad\text{for $t\not\equiv 0\pmod{p}$.}
\end{equation}

Now fix $b_0\in S^0$ such that $p\nmid\norm b_0$. Find $x_i\in
I_i$ such that $p\nmid\norm x_i/\norm I_i$; since $B$ is ramified
at $p$, the maximal order at $p$ is unique, and so $b_i :=
x_i\,b_0\,x_i^{-1}\in S^0(\Z_l) = S_i^0(\Z_p)$. We define
\[
   \wl_p(I_i,b) := \wl_{\psi,b_i}(b),
\]
to be the weight function for the triple
$(S_i^0(\Z_p),b_i,\psi)$. In practice, one can use the same
$b_0$ and $b_i$ for the definitions of both $\wl_l(I_i,\cdot)$ and
$\wl_p(I_i,\cdot)$.

Note that different choices of $\psi$ will, in general, yield
different forms $\thetaop_{-l}([I_i])$, but as long
as~\eqref{eq:psimod1} holds their coefficients will be the same
up to a constant of absolute value $1$; thus
formula~\eqref{eq:formula} will not be affected.
Moreover, given two such odd periodic functions is not difficult
to produce another \emph{periodic} function $\chi$ with the property
that the ratio of the $m$-th coefficients of the respective
theta series will be $\chi(m)$.

The case when $\psi$ is actually a character of conductor $p$
is of particular interest, since the generalized theta series
$\thetaop_{-l}([I_i])$ will be a modular form of level $4p^2$ and
character $\psi_1$, where $\psi_1(m)=\kro{-1}{m}\psi(m)$.
From a computational point of view,
however, it will always be preferable to choose a real $\psi$, whose
values will be either $0$ or $\pm 1$, and so that the coefficients of
$\thetaop_{-l}([I_i])$ will be rational integers.
Only in case $p\equiv 3\pmod{4}$ it is possible to satisfy both
requirements at the same time, by taking for $\psi$ the quadratic
character of conductor $p$.


\section{Examples}

\subsection{11A}

Let $f = f_{11A}$, the modular form of level $11$, and consider
$B=B(-1,-11)$, the quaternion algebra ramified precisely at $\infty$ and $11$.
A maximal order, and representatives for its right ideals classes, are given by
\begin{align*}
   R=I_{1}&=\<1,i,\frac{1+j}{2},\frac{i+k}{2}>
     & \text{ with $\norm I_{1} = 1$,} \\
     I_{2}&=\<2,2i,\frac{1+2i+j}{2},\frac{2+3i+k}{2}>
     & \text{ with $\norm I_{2} = 2$.} \\
\end{align*}
By computing the Brandt matrices (see \S6 of~\cite{g} for this
example), we find a vector
\[
  e_f = [I_2]-[I_1]
\]
of height $\<e_f,e_f>=5$ corresponding to $f$.
Since $L(f,1)\approx 0.25384186$, Gross's method works,
and it's easy to compute
\[
  \thetaop_1(e_f)=\thetaop_1([I_2])-\thetaop_1([I_1])
  = q^{3}-q^{4}-q^{11}-q^{12}+q^{15}+2q^{16}+O(q^{20}),
\]
as the difference of two regular theta series corresponding
to the ternary quadratic forms~\eqref{eq:11A:Q1} and~\eqref{eq:11A:Q2}.

\subsubsection{Real twists in a case of rank $0$}
Let $l=3$. One can compute $L(f,-3,1)\approx 1.6844963$,
and thus expect $\thetaop_{-3}(e_f)$ to be nonzero.
We can choose $b_0=i+k\in S^0$ with norm $12$,
and let $\psi=\chi_{11}$
be the quadratic character of conductor $11$.

Clearly we can take $x_1=1$ and $x_2=2$, so that $n_1=1$, $n_2=2$
and $b_1=b_2=i+k$. Bases for $S^0_1$ and $S^0_2$ are given by
\begin{align*}
S^0_1 &= \<2i,j,i+k>               && \text{ with $b_1 = (0,0,1)$,} \\
S^0_2 &= \<4i,2i+j,\frac{7i+k}{2}> && \text{ with $b_2 = (-\tf{3}{2},0,2)$.}
\end{align*}
The norm form in the given bases will be
\begin{align}
  \label{eq:11A:Q1}
  \norm_1{}(x_1,x_2,x_3) & =
    4 x_1^2 + 11 x_2^2 + 12 x_3^2 + 4 x_1 x_3, \\
  \label{eq:11A:Q2}
  \norm_2{}(x_1,x_2,x_3) & =
    16 x_1^2 + 15 x_2^2 + 15 x_3^2
    + 14 x_2 x_3 + 28 x_1 x_3 + 16 x_1 x_2.
\end{align}
This information is all that we need to compute $\thetaop_{-3}$.
As an example, we show how to compute $\thetaop_{-3}([I_1])$.
A simple calculation shows that
\[
  \<(x_1,x_2,x_3),b_1> = 4x_1+24x_3 \equiv x_1\pmod{3},
\]
and thus $\wl_3(I_1,\cdot)$ can be computed by
\begin{align*}
  \wl_3(I_1,(x_1,x_2,x_3)) &= \begin{cases}
    0
    & \text{if $3\nmid\norm_1{}(x_1,x_2,x_3)$,} \\
    \chi_3(x_1)
    & \text{if $x_1\not\equiv 0\pmod{3}$,} \\
    \chi_3(x_3)
    & \text{otherwise.}
  \end{cases} \\
\intertext{
Similarly, $\wl_{11}(I_1,\cdot)$ will be given by
}
  \wl_{11}(I_1,(x_1,x_2,x_3)) &= \chi_{11}(4x_1+2x_3).
\end{align*}
Hence we compute
\[
  \thetaop_{-3}([I_1]) =
    -2q^{4}+2q^{5}+2q^{9}+2q^{12}+2q^{20}+2q^{25}-2q^{37}+O(q^{48}).
\]
In a similar way one can easily get
\[
  \thetaop_{-3}([I_2]) =
    q+q^{4}-3q^{5}-3q^{12}+4q^{16}-3q^{20}+2q^{25}-6q^{36}+3q^{37}+O(q^{48}).
\]

\begin{table}
\begin{tabular}{|r|rr|r|rr|r|rr|}
\hline
$d$ & $c_{-3}({d})$ & $L(f,d,1)$& $d$ & $c_{-3}({d})$ & $L(f,d,1)$&
$d$ & $c_{-3}({d})$ & $L(f,d,1)$\\
\hline
   1 & 1 & 0.253842 & 92 & -5 & 0.661621 & 141 & -10 & 2.137734 \\
   5 & -5 & 2.838038 & 93 & 5 & 0.658054 & 152 & -10 & 2.058929 \\
   12 & -5 & 1.831946 & 97 & 5 & 0.644343 & 157 & -15 & 4.558227 \\
   37 & 5 & 1.043284 & 104 & 10 & 2.489124 & 168 & 10 & 1.958432 \\
   53 & 10 & 3.486786 & 113 & -5 & 0.596986 & 177 & 5 & 0.476998 \\
   56 & 10 & 3.392105 & 124 & -5 & 0.569892 & 181 & -15 & 4.245281 \\
   60 & -5 & 0.819271 & 133 & 10 & 2.201088 & 185 & -5 & 0.466571 \\
   69 & 15 & 6.875768 & 136 & 10 & 2.176676 & 188 & -10 & 1.851332 \\
   89 & -5 & 0.672680 & 137 & -5 & 0.542179 & & &  \\
\hline
\end{tabular}
\caption{Coefficients of $\thetaop_{-3}(e_f)$ and central values for
         $f=f_{11A}$\label{table:coeff11A}}
\end{table}
Table~\ref{table:coeff11A} shows the values of $c_{-3}(d)$ and $L(f,d,1)$,
where $0<d<200$ is a fundamental discriminant such that
$\kro{d}{11} = 1$.
The formula
\[
  L(f,d,1) = k_{-3}\,\frac{c_{-3}(d)^2}{\sqrt{d}}
\]
is satisfied, where
\[
  k_{-3} = \frac{1}{5}\cdot\frac{(f,f)}{L(f,-3,1)\sqrt{3}}
    = L(f,1)
    \approx 0.2538418608559106843377589233509...
\]
Note that when $\kro{d}{11}\neq 1$ it is
trivial that $c_{-3}(d)=L(f,d,1)=0$.


\subsection{37A}

Let $f = f_{37A}$, the modular form of level $37$ and rank $1$, and consider
$B=B(-2,-37)$, the quaternion algebra ramified precisely at $\infty$ and $37$.
A maximal order, and representatives for its right ideal classes, are given by
\begin{align*}
   R=I_1&=\<1,i,\frac{1+i+j}{2},\frac{2+3i+k}{4}>
     & \text{with $\norm I_1 = 1$,} \\
     I_2&=\<2,2i,\frac{1+3i+j}{2},\frac{6+3i+k}{4}>
     & \text{with $\norm I_2 = 2$,} \\
     I_3&=\<4,2i,\frac{3+3i+j}{2},\frac{6+i+k}{2}>
     & \text{with $\norm I_3 = 4$.} \\
\end{align*}
By computing the Brandt matrices, we find a vector
\[
   e_f=\frac{[I_3]-[I_2]}{2}
\]
of height $\<e_f,e_f>=\tf{1}{2}$ corresponding to $f$.
Since $L(f,1)=0$ we know that
$2\thetaop_1(e_f)=\thetaop_1([I_3])-\thetaop_1([I_2])=0$.
Indeed, one checks that $R_2$ and $R_3$ are conjugate, which explains
the identity $\thetaop_1([I_2])=\thetaop_1([I_3])$.

\subsubsection{Imaginary twists in a case of rank $1$}
Let $l=5$. One can compute $L(f,5,1)\approx 5.3548616$, and thus we
expect $\thetaop_5(e_f)$ to be nonzero.
We note that, by the same reason that the orders are conjugate, we
have $\thetaop_5([I_3])=-\thetaop_5([I_2])$,
except now there's an extra sign,
ultimately coming from the fact that $\kro{37}{5}=-1$.
Thus, $\thetaop_5(e_f)=\thetaop_5([I_3])$.
A basis for $S_3^0$ is given by
\[
   S_3^0 = \<4i,3i+j,\frac{3i+2j+k}{4}>,
\]
with the norm in this basis
\[
  \norm_3{}(x_1,x_2,x_3) =
    32 x_1^2 + 55 x_2^2 + 15 x_3^2
   +46 x_2 x_3 +12 x_1 x_3 +48 x_1 x_2.
\]
Choose $b_3=(0,0,1)$, with norm $15$. Then
\begin{equation}\label{eq:37A:trace}
  \<(x_1,x_2,x_3),b_3> = 12x_1+46x_2+30x_3 \equiv 2x_1+x_2\pmod{5},
\end{equation}
so that
\[
  \wl_5(I_3,(x_1,x_2,x_3)) = \begin{cases}
    0
    & \text{if $5\nmid\norm_3{}(x_1,x_2,x_3)$,} \\
    \chi_5(2x_1+x_2)
    & \text{if $2x_1+x_2\not\equiv 0\pmod{5}$,} \\
    \chi_5(x_3)
    & \text{otherwise.}
  \end{cases}
\]

\begin{table}
\begin{tabular}{|r|rr|r|rr|r|rr|}
\hline
$-d$ & $c_{5}(d)$ & $L(f,-d,1)$& $-d$ & $c_{5}(d)$ & $L(f,-d,1)$& $-d$ &
$c_{5}(d)$ & $L(f,-d,1)$\\
\hline
   -3 & 1 & 2.830621 & -95 & 0 & 0.000000 & -139 & 0 & 0.000000 \\
   -4 & 1 & 2.451389 & -104 & 0 & 0.000000 & -148 & -3 & 7.254107 \\
   -7 & -1 & 1.853076 & -107 & 0 & 0.000000 & -151 & -2 & 1.595930 \\
   -11 & 1 & 1.478243 & -111 & 1 & 0.930702 & -152 & -2 & 1.590671 \\
   -40 & 2 & 3.100790 & -115 & -6 & 16.458713 & -155 & 2 & 1.575203 \\
   -47 & -1 & 0.715144 & -120 & -2 & 1.790242 & -159 & 1 & 0.388816 \\
   -67 & 6 & 21.562911 & -123 & 3 & 3.978618 & -164 & -1 & 0.382843 \\
   -71 & 1 & 0.581853 & -127 & 1 & 0.435051 & -184 & 0 & 0.000000 \\
   -83 & -1 & 0.538150 & -132 & 3 & 3.840589 & -195 & 2 & 1.404381 \\
   -84 & -1 & 0.534937 & -136 & 4 & 6.726557 & & &  \\
\hline
\end{tabular}
\caption{Coefficients of $\thetaop_5(e_f)$ and central values for
         $f=f_{37A}$\label{table:coeff37A}}
\end{table}
Table~\ref{table:coeff37A} shows the values of $c_5(d)$ and $L(f,-d,1)$,
where $-200<-d<0$ is a fundamental discriminant such that
$\kro{-d}{37}\neq -1$.
The formula
\[
  L(f,-d,1) = k_5\,\frac{c_5(d)^2}{\sqrt{d}} \cdot
   \begin{cases}
    1 & \text{if $\kro{-d}{37} = +1$,} \\
    2 & \text{if $\kro{-d}{37} = 0$,} \\
    0 & \text{if $\kro{-d}{37} = -1$,} \\
  \end{cases}
\]
is satisfied, where
\[
  k_5 = 2\cdot\frac{(f,f)}{L(f,5,1)\sqrt{5}}
    \approx 4.902778763973580121708449663733...
\]
Note that in the case $\kro{-d}{37}=-1$ it is
trivial that $c_5(d)=L(f,-d,1)=0$.

\subsubsection{Real twists in a case of rank $1$}
Let $l=3$, since $L(f,-3,1)\approx 2.9934586$.
Keep $b_3$ as above, and let $\psi$ be
the odd periodic function modulo $37$ such that
\[
   \psi(x) = \begin{cases}
     +1 & \text{ if $\phantom{0}1\leq x\leq 18$,}\\
     -1 & \text{ if $19\leq x\leq 36$.}
   \end{cases}
\]
Using again~\eqref{eq:37A:trace}, we have that
\begin{align*}
  \wl_3(I_1,(x_1,x_2,x_3)) &= \begin{cases}
    0
    & \text{if $3\nmid\norm_1{}(x_1,x_2,x_3)$,} \\
    \chi_3(x_2)
    & \text{if $x_2\not\equiv 0\pmod{3}$,} \\
    \chi_3(3)
    & \text{otherwise.}
  \end{cases} \\
\intertext{
and $\wl_{11}(I_1,\cdot)$ will be given by
}
  \wl_{37}(I_1,(x_1,x_2,x_3)) &= \psi(12x_1+9x_2+30x_3).
\end{align*}

\begin{table}
\begin{tabular}{|r|rr|r|rr|r|rr|}
\hline
$d$ & $c_{-3}(d)$ & $L(f,d,1)$&
$d$ & $c_{-3}(d)$ & $L(f,d,1)$&
$d$ & $c_{-3}(d)$ & $L(f,d,1)$\\
\hline
   5 & 1 & 5.354862 & 76 & 1 & 1.373493 & 133 & -1 & 1.038263 \\
   8 & 1 & 4.233390 & 88 & 1 & 1.276415 & 140 & -3 & 9.107764 \\
   13 & -1 & 3.320944 & 89 & -1 & 1.269224 & 156 & -1 & 0.958674 \\
   17 & 1 & 2.904081 & 92 & 2 & 4.993434 & 161 & -2 & 3.774681 \\
   24 & -1 & 2.444149 & 93 & 2 & 4.966515 & 165 & 1 & 0.932162 \\
   29 & 2 & 8.893941 & 97 & 0 & 0.000000 & 168 & -1 & 0.923801 \\
   56 & -1 & 1.600071 & 105 & 1 & 1.168527 & 172 & 1 & 0.912996 \\
   57 & 1 & 1.585973 & 109 & -1 & 1.146885 & 177 & 0 & 0.000000 \\
   60 & -1 & 1.545815 & 113 & 0 & 0.000000 & 193 & -1 & 0.861895 \\
   61 & 0 & 0.000000 & 124 & 0 & 0.000000 & & &  \\
   69 & 0 & 0.000000 & 129 & 1 & 1.054237 & & &  \\
\hline
\end{tabular}
\caption{Coefficients of $\thetaop_{-3}(e_f)$ and central values for
         $f=f_{37A}$\label{table:coeff37Are}}
\end{table}
Table~\ref{table:coeff37Are} shows the values of $c_{-3}(d)$ and $L(f,d,1)$,
where $0<d<200$ is a fundamental discriminant such that
$\kro{d}{37}=-1$.
The formula
\[
  L(f,d,1) = k_{-3}\,\frac{c_{-3}(d)^2}{\sqrt{d}}
\]
is satisfied, where
\[
  k_{-3} = 2\cdot\frac{(f,f)}{L(f,-3,1)\sqrt{3}}
    \approx 11.97383458492783851932803991781...
\]
Note that in the case $\kro{d}{37}\neq-1$ it is
trivial that $c_{-3}(d)=L(f,d,1)=0$.

\subsection{43A}

Let $f = f_{43A}$, the modular form of level $43$ and rank $1$.
Let $B=B(-1,-43)$,
the quaternion algebra ramified precisely at $\infty$ and $43$.
A maximal order, and representatives for its right ideals classes, are
given by
\begin{align*}
   R=I_1&=\<1,i,\frac{1+j}{2},\frac{i+k}{2}>
     & \text{with $\norm I_1 = 1$,} \\
     I_2&=\<2,2i,\frac{1+2i+j}{2},\frac{2+3i+k}{2}>
     & \text{with $\norm I_2 = 2$,} \\
     I_3&=\<3,3i,\frac{1+2i+j}{2},\frac{2+5i+k}{2}>
     & \text{with $\norm I_3 = 3$,} \\
     I_4&=\<3,3i,\frac{1+4i+j}{2},\frac{4+5i+k}{2}>
     & \text{with $\norm I_4 = 3$.} \\
\end{align*}
By computing the Brandt matrices, we find a vector
\[
  e_f = \frac{[I_4]-[I_3]}{2}
\]
of height $\<e_f,e_f>=\tf{1}{2}$ corresponding to $f$.

\subsubsection{Imaginary twists in a case of rank $1$}
We can use $l=5$, since $L(f,5,1)\approx 4.8913446$
is nonzero; again, we find $\thetaop_5(e_f)=\thetaop_5([I_4])$.
\begin{table}
\begin{tabular}{|r|rr|r|rr|r|rr|}
\hline
$-d$ & $c_{5}(d)$ & $L(f,-d,1)$& $-d$ & $c_{5}(d)$ & $L(f,-d,1)$& $-d$ &
$c_{5}(d)$ & $L(f,-d,1)$\\
\hline
   -3 & 1 & 3.148135 & -91 & -1 & 0.571601 & -151 & -1 & 0.443737 \\
   -7 & 1 & 2.060938 & -104 & 1 & 0.534684 & -155 & -1 & 0.437974 \\
   -8 & -1 & 1.927831 & -115 & -3 & 4.576227 & -159 & 1 & 0.432430 \\
   -19 & 2 & 5.003768 & -116 & -1 & 0.506273 & -163 & 7 & 20.927447 \\
   -20 & -1 & 1.219267 & -119 & -1 & 0.499851 & -168 & -2 & 1.682749 \\
   -39 & -1 & 0.873136 & -120 & 0 & 0.000000 & -179 & -1 & 0.407556 \\
   -43 & 2 & 6.652268 & -123 & -5 & 12.291402 & -184 & -3 & 3.617825 \\
   -51 & 1 & 0.763535 & -131 & 0 & 0.000000 & -191 & 0 & 0.000000 \\
   -55 & 1 & 0.735246 & -132 & 3 & 4.271393 & -199 & 0 & 0.000000 \\
   -71 & 0 & 0.000000 & -136 & -1 & 0.467568 & & &  \\
   -88 & 3 & 5.231366 & -148 & -4 & 7.171386 & & &  \\
\hline
\end{tabular}
\caption{Coefficients of $\thetaop_5(e_f)$ and central values for
         $f=f_{43A}$\label{table:coeff43A}}
\end{table}
Table~\ref{table:coeff43A} shows the values of $c_5(d)$ and $L(f,-d,1)$,
where $-200<-d<0$ is a fundamental discriminant such that
$\kro{-d}{43}\neq -1$.
The formula
\[
  L(f,-d,1) = k_5\,\frac{c_5(d)^2}{\sqrt{d}} \cdot
   \begin{cases}
    1 & \text{if $\kro{-d}{43} = +1$,} \\
    2 & \text{if $\kro{-d}{43} = 0$,} \\
    0 & \text{if $\kro{-d}{43} = -1$,} \\
  \end{cases}
\]
is satisfied, where
\[
  k_5 = 2\cdot\frac{(f,f)}{L(f,5,1)\sqrt{5}}
    \approx 5.452729672681734385570722785283...
\]
Note that in the case $\kro{-d}{43}=-1$ it is
trivial that $c_5(d)=L(f,-d,1)=0$.

\subsubsection{Real twists in a case of rank $1$}
We can use $l=3$, since $L(f,-3,1)\approx 3.1481349$,
and let $\psi=\chi_{43}$ be the quadratic character of conductor $43$.
\begin{table}
\begin{tabular}{|r|rr|r|rr|r|rr|}
\hline
$d$ & $c_{-3}(d)$ & $L(f,d,1)$&
$d$ & $c_{-3}(d)$ & $L(f,d,1)$&
$d$ & $c_{-3}(d)$ & $L(f,d,1)$\\
\hline
   5 & 1 & 4.891345 & 76 & 0 & 0.000000 & 137 & 2 & 3.737773 \\
   8 & -1 & 3.866947 & 77 & -3 & 11.217870 & 141 & -2 & 3.684374 \\
   12 & 1 & 3.157349 & 85 & 1 & 1.186325 & 149 & 0 & 0.000000 \\
   28 & -1 & 2.066970 & 88 & -1 & 1.165929 & 156 & 1 & 0.875691 \\
   29 & -1 & 2.031020 & 89 & 1 & 1.159360 & 157 & 2 & 3.491592 \\
   33 & -1 & 1.903953 & 93 & 3 & 10.207380 & 161 & -1 & 0.861986 \\
   37 & 2 & 7.192376 & 104 & 1 & 1.072498 & 168 & -2 & 3.375348 \\
   61 & 1 & 1.400388 & 105 & 0 & 0.000000 & 177 & -2 & 3.288415 \\
   65 & -1 & 1.356615 & 113 & -2 & 4.115608 & 184 & 1 & 0.806314 \\
   69 & -1 & 1.316706 & 120 & 0 & 0.000000 & & &  \\
   73 & 1 & 1.280123 & 136 & 1 & 0.937873 & & &  \\
\hline
\end{tabular}
\caption{Coefficients of $\thetaop_{-3}(e_f)$ and central values for
         $f=f_{43A}$\label{table:coeff43Are}}
\end{table}
Table~\ref{table:coeff43Are} shows the values of $c_{-3}(d)$ and $L(f,d,1)$,
where $0<d<200$ is a fundamental discriminant such that
$\kro{d}{43}=-1$.
The formula
\[
  L(f,d,1) = k_{-3}\,\frac{c_{-3}(d)^2}{\sqrt{d}}
\]
is satisfied, where
\[
  k_{-3} = 2\cdot\frac{(f,f)}{L(f,-3,1)\sqrt{3}}
    \approx 10.937379059935167648758735438779...
\]
Note that in the case $\kro{d}{43}\neq-1$ it is
trivial that $c_{-3}(d)=L(f,d,1)=0$.

\subsection{389A}

Let $f = f_{389A}$, the modular form of level $389$ and rank $2$.
Let $B=B(-2,-389)$,
the quaternion algebra ramified precisely at $\infty$ and $389$.
A maximal order, with $33$ ideal classes, is given by
\[
   R=\<1,i,\frac{1+i+j}{2},\frac{2+3i+k}{4}>.
\]
There is a vector $e_f$
of height $\<e_f,e_f>=\tf{5}{2}$
corresponding to $f$.

\subsubsection{Imaginary twists in a case of rank $2$}
We can use $l=5$, since $L(f,5,1)\approx 8.9092552$.
We have omitted the $33$ ideal classes; however,
the computation of $\thetaop_l(e_f)$ involves only $14$ distinct theta
series. In table~\ref{table:qf389A} we give the value of $e_f$ and the
coefficients of the norm form $\norm_i{}$ and of $b_i$ on chosen bases of $S_i^0$.

\begin{table}
\begin{tabular}{|c|r|r@{, }r@{, }r@{, }r@{, }r@{, }r|r@{, }r@{, }r|}
  \hline
  $i$ & $a_i$ &
  \multicolumn{6}{c|}{$\norm_i{}$} &
  \multicolumn{3}{c|}{$b_i$} \\
  \hline
   1 &  1/2 & 15&107&416&-100& -8&-14  &  2&4&0 \\
   2 & -1/2 & 15&104&415& 104&  2&  4  &  0&4&1 \\
   3 & -1/2 & 23&136&203&  68&  2&  8  &  2&1&4 \\
   4 &  1/2 & 23& 72&407&  72& 10& 20  &  1&1&0 \\
   5 & -1/2 & 31& 51&407& -46&-26&-10  &  1&2&0 \\
   6 &  1/2 & 31&103&204&  56& 20& 18  &  2&0&3 \\
   7 &  1/2 & 39&128&160&-116& -8&-36  &  1&1&4 \\
   8 & -1/2 & 39& 40&399&  40&  2&  4  &  1&0&1 \\
   9 &  1/2 & 40& 47&399&  18& 40& 36  &  4&3&0 \\
  10 & -1/2 & 47&107&135&  42& 22& 38  &  4&3&1 \\
  11 & -1/2 & 56& 84&139&  56&  4& 12  &  3&1&4 \\
  12 &  1/2 & 56& 92&151&  76& 52& 44  &  4&2&3 \\
  13 &  1/2 & 71& 83&132& -16&-12&-70  &  2&3&4 \\
  14 & -1/2 & 71&103&124& -36&-64&-66  &  4&0&2 \\
  \hline
\end{tabular}
\caption{Coefficients of the ternary forms
         and of $b_i$\label{table:qf389A}}
\end{table}
Each row in the table allows one to compute an individual theta series
\[
  h_i(z) := \frac{1}{2}\sum_{b\in\Z^3} w_5(I_i,b) q^{\norm_i{}(b)/5}.
\]
The ternary form corresponding to a sextuple
$(A_1,A_2,A_3,A_{23},A_{13},A_{12})$ is
\[
  \norm_i{}(x_1,x_2,x_3) =
    A_1 x_1^2 + A_2 x_2^2 + A_3 x_3^2
    + A_{23} x_2 x_3 + A_{13} x_1 x_3 + A_{12} x_1 x_2,
\]
and $\wl_5(I_i,\cdot)$ is the weight function of the pair
$(\Z^3,b_i)$. As an example, we show how to compute $h_1(z)$.
First, we have
\[
  \norm_1{}(x_1,x_2,x_3) =
    15 x_1^2 + 107 x_2^2 + 416 x_3^2
    - 100 x_2 x_3 - 8 x_1 x_3 -14 x_1 x_2.
\]
A simple calculation shows that
\[
  \<(x_1,x_2,x_3),(2,4,0)> \equiv 4x_1+3x_2+4x_3 \pmod{5}.
\]
Thus, $\wl_5$ can be computed as
\[
  \wl_5(I_1,(x_1,x_2,x_3)) = \begin{cases}
    0
    & \text{if $5\nmid\norm_1{}(x_1,x_2,x_3)$,} \\
    \chi_5(4x_1+3x_2+4x_3)
    & \text{if $\not\equiv 0\pmod{5}$,} \\
    \chi_5(x_2)
    & \text{otherwise,}
  \end{cases}
\]
%
and we have
\[
  h_1(z)=q^{3}-q^{12}-q^{27}+q^{39}+q^{40}+q^{48}-q^{83}-2q^{92}+O(q^{100}).
\]
Finally, we combine all of the theta series in
\[
  \thetaop_5(e_f) = \sum_{i=1}^{14} a_i h_i(z)
\]

\begin{table}
\begin{tabular}{|r|rr|r|rr|r|rr|}
\hline
$-d$ & $c_{5}(d)$ & $L(f,-d,1)$& $-d$ & $c_{5}(d)$ & $L(f,-d,1)$& $-d$ &
$c_{5}(d)$ & $L(f,-d,1)$\\
\hline
    -3 &  1 &  4.553533 &  -83 & -1 &  0.865705 & -139 & -1 &  0.668962 \\
    -8 & -1 &  2.788458 &  -84 &  1 &  0.860537 & -148 &  6 & 23.338921 \\
   -15 & -1 &  2.036402 &  -88 & -4 & 13.452028 & -151 &  2 &  2.567324 \\
   -23 &  1 &  1.644543 & -103 &  0 &  0.000000 & -152 & -1 &  0.639716 \\
   -31 &  1 &  1.416538 & -104 & -1 &  0.773379 & -155 &  3 &  5.701456 \\
   -39 &  1 &  1.262923 & -107 &  0 &  0.000000 & -163 &  8 & 39.536232 \\
   -40 &  1 &  1.247036 & -115 & -1 &  0.735462 & -167 & -1 &  0.610311 \\
   -43 & -3 & 10.824738 & -116 & -2 &  2.929140 & -191 &  1 &  0.570680 \\
   -47 &  0 &  0.000000 & -123 &  3 &  6.400282 & -195 &  1 &  0.564796 \\
   -51 & -2 &  4.417576 & -131 &  1 &  0.689086 & -199 & -1 &  0.559091 \\
   -56 &  1 &  1.053938 & -132 & -2 &  2.745884 & & & \\ 
   -71 &  1 &  0.936009 & -136 & -2 &  2.705202 & & & \\ 
\hline
\end{tabular}
\caption{Coefficients of $\thetaop_5(e_f)$ and central values for
         $f=f_{389A}$\label{table:coeff389A}}
\end{table}
Table~\ref{table:coeff389A} shows the values of $c_5(d)$ and $L(f,-d,1)$,
where $0<-d<200$ is a fundamental discriminant such that
$\kro{-d}{389}\neq +1$.
The formula
\[
  L(f,-d,1) = k_5\,\frac{c_5(d)^2}{\sqrt{d}} \cdot
  \begin{cases}
    1 & \text{if $\kro{-d}{389} = -1$,} \\
    2 & \text{if $\kro{-d}{389} = 0$,} \\
    0 & \text{if $\kro{-d}{389} = +1$,} \\
  \end{cases}
\]
is satisfied, where
\[
  k_5 = \frac{2}{5}\cdot\frac{(f,f)}{L(f,5,1) \sqrt{5}}
    \approx 7.886950806206592817689630792605...
\]
Note that when $\kro{-d}{389}=+1$ it is
trivial that $c_5(d)=L(f,-d,1)=0$.




\begin{thebibliography}{B-F-H}
\bibitem[BM]{bm}Baruch E.M., Mao Z., {\em Central values of automorphic
$L$-functions}, preprint.
\bibitem[BFH]{bfh}Bump, D., Friedberg,
S., Hoffstein, J., {\em Nonvanishing theorems for $L$-functions of
modular forms and their derivatives}, Invent. Math. 102 (1990), p.
543-618.
\bibitem[G]{g} Gross, B., {\em Heights and the special
values of $L-$series},
 Canadian Math. Soc. Conf. Proceedings, volume 7, (1987) p.
 115-187.
\bibitem[P]{p}Pizer A., {\em An algorithm for computing modular
forms on $\Gamma_0(N)$}, J. Algebra 64 (1980), p. 340-390.
\bibitem[W]{wa1}Waldspurger J-L., {\em Sur les coefficients de Fourier
des formes modulaires de poids demi-entier}, J. Math. pures Appl.
60 (1981), p. 375-484.
\end{thebibliography}
\end{document}